\documentclass[11pt]{amsart}
\usepackage[a4paper, total={6in, 8in}]{geometry}

\usepackage[pagebackref]{hyperref}

\usepackage{amsmath,amsfonts,amssymb}
\usepackage{graphicx} 
\usepackage{tikz}\usetikzlibrary{matrix, cd, arrows}
\usepackage{tikz-cd}
\newcommand{\At}{\mathcal{A}t}
\newcommand{\Ato}{\mathcal{A}t^0}
\newcommand{\End}{\operatorname{End}}
\newcommand{\calD}{\mathcal{D}}
\newcommand{\calC}{\mathcal{C}}
\newcommand{\calM}{\mathcal{M}}
\newcommand{\calE}{\mathcal{E}}
\newcommand{\Id}{\operatorname{Id}}
\newcommand{\SL}{\operatorname{SL}}
\newcommand{\calL}{\mathcal{L}}
\newcommand{\dol}{\operatorname{Dol}}
\newcommand{\dr}{\operatorname{dR}}
\newcommand{\PP}{\mathbb{P}}
\newcommand{\twistor}{\operatorname{tw}}
\newcommand{\Gm}{\mathbb{G}_m}
\newcommand{\hod}{\operatorname{Hod}}
\newcommand{\AAA}{\mathbb{A}}
\newcommand{\calO}{\mathcal{O}}
\newcommand{\CC}{\mathbb{C}}
\newcommand{\s}{\operatorname{s}}
\newcommand{\sss}{\operatorname{ss}}

\theoremstyle{plain}
\newtheorem{theorem}{Theorem}
\newtheorem{proposition}[theorem]{Proposition}

\theoremstyle{definition}

\theoremstyle{remark}

\title{The Determinant of Cohomology and Moduli of $\lambda$-Connections}
\author{Johan Martens}
\address{Johan Martens\\ School of Mathematics and Maxwell Institute\\ The University of Edinburgh\\ Peter Guthrie Tait Road\\ Edinburgh EH9 3FD\\ United Kingdom}
\email{johan.martens@ed.ac.uk}
\date{\today}
\begin{document}
\begin{abstract}
    We exhibit how the Hodge-Deligne moduli space of $\lambda$-connections over a smooth projective curve, for stable bundles with fixed determinant, can be understood as the dual of the Atiyah algebroid of the determinant of cohomology line bundle.
\end{abstract}
\maketitle
\section{Introduction}
When studying various moduli spaces of bundles with attributes on a smooth projective curve, one finds that many of these moduli spaces are (at least on a suitable stable locus) vector bundles over a more basic moduli space, such as the moduli space of stable bundles.  The most prominent case is that of (ordinary) Higgs bundles \cite{Hitchin1,Hitchin2}, whose moduli space, restricted to those Higgs bundles with a stable underlying vector bundle, is the total space of the cotangent bundle of the moduli space of stable bundles, and in particular is complex symplectic.  

\medskip

A variation of this can arise for moduli spaces that are complex Poisson, such as moduli spaces of weakly parabolic \cite{Logares.Martens} or parahoric \cite{MR4705763} Higgs bundles, or Higgs bundles over stacky curves \cite{KydonakisSunZhao+2024+163+182}.  In these cases, the moduli spaces (for underlying stable bundles) can be seen to be the dual of the Atiyah algebroid for a principal bundle over the moduli space of parabolic / parahoric / stacky bundles.  Duals of Atiyah algebroids always have a (linear) Poisson structure, generalizing the Kirrilov--Kostant--Souriau Poisson structure on the dual of a Lie algebra.

\medskip

In this note, we exhibit another instance of this latter phenomenon, by realising the Hodge-Deligne moduli space of $\lambda$-connections as the dual to the Atiyah algebroid of the determinant of cohomology line bundle over the moduli space of stable bundles.
\section{Atiyah algebroids of vector bundles}
Let $E\rightarrow M$ be a vector bundle over a non-singular variety.  The Atiyah algebroid of $E$,  $\At(E)$, is the sheaf of first-order differential operators on $E$ that have diagonal symbol $\sigma$, i.e. the middle term in the top sequence (the Atiyah sequence) of 
$$
\begin{tikzcd}
0\ar[r]&\End(E)\ar[r]\ar[d,equal]& \At(E)\ar[r,"\sigma"]\ar[d]& T_M\ar[r]\ar[d, "-\otimes \Id"]& 0\\
0\ar[r]&\End(E)\ar[r]& \calD(E)\ar[r,"\sigma"]& T_M\otimes\End(E)\ar[r]& 0.\\
\end{tikzcd}$$
A connection on $E$ is a splitting of the Atiyah sequence \cite{MR0086359}.  There is a natural morphism $\At(E)\rightarrow \At(\det E))$, and if the determinant of $E$ is trivial (so that its Atiyah algebroid canonically splits) we can define the traceless Atiyah algebroid $\Ato(E)$ as the kernel of the trace morphism via 
$$
\begin{tikzcd}
0\ar[r]&\End^0(E)\ar[r]\ar[d,hook]& \Ato(E)\ar[r]\ar[d]& T_M\ar[r]\ar[d, equal ]& 0\\
0\ar[r]&\End(E)\ar[r]\ar[d,"\operatorname{tr}"]& \At(E)\ar[r]\ar[d,"\operatorname{tr}+\sigma"]& T_M\ar[r]\ar[d, equal ]& 0\\
0\ar[r]& \calO_M\ar[r]& \At(\det E)=\calO\oplus T_M\ar[r] & T_M\ar[r] & 0.
\end{tikzcd}
$$
As for any Lie algebroids, the total spaces of the duals of $\Ato(E)$ and $\At(E)$ (considered as vector bundles) naturally come equipped with Poisson structures.  If we further have a smooth morphism $M\rightarrow S$, there are also natural relative versions $\At_{M/S}(E)$ and $\At^0_{M/S}(E)$, whose symbols live in the relative tangent sheaf $T_{M/S}$.

\section{Atiyah algebroid of the determinant-of-cohomology}
Let $\calC$ now be a smooth projective curve.  We can consider the moduli space of stable rank $r$ bundles with trivial determinant on $\calC$, which we will denote as $\calM_{\SL_r}^{\s}$.  This is a smooth, quasiprojective variety.  For simplicity of notation, we will pretend a universal bundle $\calE$ exists on $\calM_{\SL_r}^{\s}\times \calC$ -- this bundle does not descend from the moduli stack to the moduli space, but all the associated bundles we will use do.

\medskip

With the projection $\pi:\mathcal{M}_{\SL_r}^{\s}\times \calC\rightarrow \mathcal{M}$, we can consider the determinant of cohomology line bundle \cite{MR0437541,MR0783704}, $$\calL_{\det}=\det \pi_*\calE\otimes \det R^1\pi_*(\calE)^*.$$
This line bundle generates the Picard group of $\calM_{\SL_r}^{\s}$.  In complex differential geometry, the Chern class of this line bundle is given via the Narasimhan-Atiyah-Bott Kaehler form, which is the curvature of a connection on $\calL_{\det}$.  

\medskip

The same ingredients used to write down this differential form (the Killing form for $\SL_r$, and Serre duality / integration over $\calC$) also allow one to write down the Atiyah sequence of $\calL_{\det}$, the extension class of which is the Chern class.  In full generality, for arbitrary families of bundles on $\calC$ without fixing the determinant, this is handled by the theory of \emph{trace complexes} of Beilinson and Schechtman \cite{MR0962493, MR1977584, MR1754519}.  

\medskip

Trace complexes are quite complicated, but this theory simplifies significantly though for stable bundles with trivial determinant (see \cite[Remark 9.2]{MR1363059}, \cite[Lemma 4.11]{MR2106127}, \cite[Theorem 4.4.1]{MR4556938}): 
\begin{theorem}
There exists a isomorphism of short exact sequences over $\mathcal{M}_{\SL_R}^{\s}$
\begin{equation}\label{tracesimple}\begin{tikzcd}
0\ar[r]& R^1\pi_* T_{\calC}^*\ar[r]\ar[d, "\cong"]& R^1\pi_* \left(\Ato_{\calM\times\calC/\calM}(\calE)\right)^*\ar[r]\ar[d,"\cong"]&  R^1\pi_*\End^0(\calE)^*\ar[r]\ar[d,"\cong"]& 0\\
0\ar[r]& \calO_{\mathcal{M}^{\s}_{\SL_r}}\ar[r]& \At(\calL_{\det})\ar[r]& T_{\mathcal{M}^{\s}_{\SL_r}}\ar[r]& 0,
\end{tikzcd}\end{equation} where the left vertical isomorphism is given by relative Serre duality, and the right vertical isomorphism uses the Killing form to identify $\End^0(\mathcal{E})^*\cong \End^0(\mathcal{E})$.
\end{theorem}

\section{The Hodge-Deligne moduli space of $\lambda$-connections}
The moduli space $\calM_{\dol}^{\sss}$ of semi-stable, rank $r$ Higgs bundles with fixed determinant on $\calC$ has a hyper-K\"ahler structure, and if we denote the complex structure of $\calM_{\dol}$ as $I$, a complex structure $J$ has a modular interpretation $\calM_{\dr}$ as the moduli space $\calM_{\dr}$ of flat $\SL_r$ connections (hence there is a homomorphism $\calM^{\sss}_{\dol}\cong \calM_{\dr}$, part of non-abelian Hodge theory).

\medskip

One can therefore look at the \emph{twistor space} of $\calM^{\sss}_{\dol}$, which is a complex analytic Poisson variety $\calM_{\twistor}\rightarrow \PP^1$, whose fibres are complex symplectic, giving the various complex structures of the hyper-K\"ahler structure \cite[\S 3.F]{MR0877637}.

\medskip

Following a suggestion by Deligne, Simpson described the twistor space $\calM_{\twistor}$ as follows \cite[\S 4]{MR1492538}: in this case, there is a $\Gm$ action on $\calM_{\twistor}$, and one can understand $\calM_{\twistor}$ by glueing the Hodge-Deligne moduli space $\calM_{\hod,\calC}$ to $\calM_{\hod,\overline{\calC}}$, where $\overline{\calC}$ is the curve obtained by taking the opposite complex structure on $\calC$.  The Hodge-Deligne moduli space, which naturally comes with a morphism $\calM_{\hod,\calC}\rightarrow \AAA^1$, arises from the moduli problem of $\lambda$-connections \cite{MR1307297}, which are triples $(E,\lambda,\nabla)$, where $E$ is a bundle, $\lambda\in \mathbb{C}$ a scalar, and $\nabla$ a morphism $E\rightarrow E\otimes T^*$ satisfying a scaled Leibniz rule:
\begin{equation}\label{scaledleib}\nabla(fs)=f\nabla(s)+\lambda df \otimes s,\end{equation} for regular functions $f$ and sections $s$ of $E$.

\medskip

When $\lambda=0$, a $\lambda$-connection is just a Higgs field, and so the $0$-fibre of $\calM_{\hod}$ is just $\calM_{\dol}$. When $\lambda=1$, we get ordinary connections, so the fibre over $1$ is just $\calM_{\dr}$, and by using the $\Gm$-action all the other fibres for $\lambda\neq 0$ are also isomorphic to this.

\medskip

In the particular case when the base is a projective curve, one can also understand $\lambda$-connections as pairs $(E,D)$, where $E$ is a bundle, and $D\in H^0(\calC,\At(E)\otimes K_{\calC})$ is a section of the Atiyah algebroid of $E$ times the canonical bundle of the curve.  By applying the symbol map, we have a morphism 
$$\begin{tikzcd} H^0(\At(E)\otimes K_{\calC})\ar[r,"\sigma\otimes\Id"] & H^0(T_{\calC}\otimes K_{\calC})\cong H^0(\calO_{\calC})=\CC\ :\   D\ar[r,mapsto]& \lambda .\end{tikzcd}$$  As the symbol of sections of $\At(E)$ is diagonal, one recovers (\ref{scaledleib}) for $D$ mapping to a given $\lambda$.
There are of course variants for fixed determinant as well.

\section{The Hodge-Deligne Moduli Space over the Stable Locus}
Let us now focus on the fixed-determinant case, and further restrict all moduli spaces to the locus where the underlying bundle is stable.  In this case we have natural morphisms $$\calM_{\dol}^{\s}\cong T^*\calM_{\SL_r}^{\s}\rightarrow \calM_{\SL_r}^{\s}\qquad \text{and}\qquad \calM_{\hod}^{\s}\rightarrow \calM^{\s}_{\SL_r},$$
which can be interpreted as vector bundles over $\calM^{\s}_{\SL_r}$.  Together with the symbol map sending a $\lambda$-connection to $\lambda$, we can interpret this as a short exact sequence 
\begin{equation}\label{dualshort}
\begin{tikzcd}
0\ar[r] 
& \calM_{\dol}^{\s}\cong T^*\calM_{\SL_r}^{\s}\ar[r] 
& \calM_{\hod}^{\s}\cong \pi_*\left(\Ato(\calE)\otimes K_{\calC}\right)\ar[r] 
& \calO_{\calM^{\s}_{\SL_r}}\ar[r]
& 0.
\end{tikzcd}\end{equation}
This has the appearance of the dual of the Atiyah sequence of a line bundle, and by looking at the middle terms of (\ref{tracesimple}) and (\ref{dualshort}) and applying relative Serre duality, we obtain
\begin{proposition} Over the locus of stable underlying vector bundles, the Hodge-Deligne space $\calM_{\hod}^{\s}$ for fixed determinants is dual to the total space of the Atiyah algebroid $\Ato(\calL_{\det})$:
$$\begin{tikzcd}
\calM_{\hod}^{\s} \ar[rr, "\cong" ]\ar[dr]& & \Ato(\calL_{\det})^*\ar[dl]\\
& \calM_{\SL_r}^{\s}.
\end{tikzcd}$$
This isomorphism moreover preserves the Poisson structures on the total spaces of these bundles.
\end{proposition}
\section*{Acknowledgements}
The author would like to thank the Isaac Newton Institute for Mathematical Sciences, Cambridge (funded by EPSRC grant EP/R014604/1), for support and hospitality during the programme \emph{New equivariant methods in algebraic and differential geometry}, where the work on this note was undertaken.  The observation made in this note arose during this programme in discussions with Jochen Heinloth, and the author would like to thank him for a generous exchange of ideas.  The author also wants to thank Lucien Hennecart and Dimitry Wyss for useful discussions.

\def\cprime{$'$} \def\cftil#1{\ifmmode\setbox7\hbox{$\accent"5E#1$}\else \setbox7\hbox{\accent"5E#1}\penalty 10000\relax\fi\raise 1\ht7 \hbox{\lower1.15ex\hbox to 1\wd7{\hss\accent"7E\hss}}\penalty 10000 \hskip-1\wd7\penalty 10000\box7} \def\cftil#1{\ifmmode\setbox7\hbox{$\accent"5E#1$}\else \setbox7\hbox{\accent"5E#1}\penalty 10000\relax\fi\raise 1\ht7 \hbox{\lower1.15ex\hbox to 1\wd7{\hss\accent"7E\hss}}\penalty 10000 \hskip-1\wd7\penalty 10000\box7}


\begin{thebibliography}{PMBB23}
\expandafter\ifx\csname url\endcsname\relax
  \def\url#1{\texttt{#1}}\fi
\expandafter\ifx\csname doi\endcsname\relax
  \def\doi#1{\burlalt{doi:#1}{http://dx.doi.org/#1}}\fi
\expandafter\ifx\csname urlprefix\endcsname\relax\def\urlprefix{URL }\fi
\expandafter\ifx\csname href\endcsname\relax
  \def\href#1#2{#2}\fi
\expandafter\ifx\csname burlalt\endcsname\relax
  \def\burlalt#1#2{\href{#2}{#1}}\fi

\bibitem[Ati57]{MR0086359}
M.~F. Atiyah.
\newblock Complex analytic connections in fibre bundles.
\newblock {\em Trans. Amer. Math. Soc.}, 85:181--207, 1957, \doi{10.2307/1992969}.

\bibitem[BE02]{MR1977584}
S.~Bloch and H.~Esnault.
\newblock Relative algebraic differential characters.
\newblock In {\em Motives, polylogarithms and {H}odge theory, {P}art {I} ({I}rvine, {CA}, 1998)}, volume 3, I of {\em Int. Press Lect. Ser.}, pages 47--73. Int. Press, Somerville, MA, 2002.
\newblock \burlalt{\tt [arxiv:math/9912015]}{http://arxiv.org/abs/math/9912015}.

\bibitem[BS88]{MR0962493}
A.~A. Be\u{\i}linson and V.~V. Schechtman.
\newblock Determinant bundles and {V}irasoro algebras.
\newblock {\em Comm. Math. Phys.}, 118(4):651--701, 1988.
\newblock \urlprefix\url{http://projecteuclid.org/euclid.cmp/1104162170}.

\bibitem[ET00]{MR1754519}
H.~Esnault and I.-H. Tsai.
\newblock Determinant bundle in a family of curves, after {A}. {B}eilinson and {V}. {S}chechtman.
\newblock {\em Comm. Math. Phys.}, 211(2):359--363, 2000, \doi{10.1007/s002200050816}.
\newblock \burlalt{\tt [arxiv:math/9910053]}{http://arxiv.org/abs/math/9910053}.

\bibitem[Gin95]{MR1363059}
V.~Ginzburg.
\newblock Resolution of diagonals and moduli spaces.
\newblock In {\em The moduli space of curves ({T}exel {I}sland, 1994)}, volume 129 of {\em Progr. Math.}, pages 231--266. Birkh\"auser Boston, Boston, MA, 1995, \doi{10.1007/978-1-4612-4264-2\_9}.
\newblock \burlalt{\tt [arxiv:hep-th/9410055]}{http://arxiv.org/abs/hep-th/9410055}.

\bibitem[Hit87a]{Hitchin1}
N.~J. Hitchin.
\newblock The self-duality equations on a {R}iemann surface.
\newblock {\em Proc. London Math. Soc. (3)}, 55(1):59--126, 1987, \doi{10.1112/plms/s3-55.1.59}.

\bibitem[Hit87b]{Hitchin2}
N.~Hitchin.
\newblock Stable bundles and integrable systems.
\newblock {\em Duke Math. J.}, 54(1):91--114, 1987, \doi{10.1215/S0012-7094-87-05408-1}.

\bibitem[HKLR87]{MR0877637}
N.~J. Hitchin, A.~Karlhede, U.~Lindstr\"om, and M.~Ro\v{c}ek.
\newblock Hyper-{K}\"ahler metrics and supersymmetry.
\newblock {\em Comm. Math. Phys.}, 108(4):535--589, 1987.
\newblock \urlprefix\url{http://projecteuclid.org/euclid.cmp/1104116624}.

\bibitem[KM76]{MR0437541}
F.~F. Knudsen and D.~Mumford.
\newblock The projectivity of the moduli space of stable curves. {I}. {P}reliminaries on ``det'' and ``{D}iv''.
\newblock {\em Math. Scand.}, 39(1):19--55, 1976, \doi{10.7146/math.scand.a-11642}.

\bibitem[KSZ24a]{MR4705763}
G.~Kydonakis, H.~Sun, and L.~Zhao.
\newblock Logahoric {H}iggs torsors for a complex reductive group.
\newblock {\em Math. Ann.}, 388(3):3183--3228, 2024, \doi{10.1007/s00208-023-02605-x}.
\newblock \burlalt{\tt [arxiv:2107.01977]}{http://arxiv.org/abs/2107.01977}.

\bibitem[KSZ24b]{KydonakisSunZhao+2024+163+182}
G.~Kydonakis, H.~Sun, and L.~Zhao.
\newblock Poisson structures on moduli spaces of {H}iggs bundles over stacky curves.
\newblock {\em Advances in Geometry}, 24(2):163--182, 2024, \doi{10.1515/advgeom-2024-0004}.
\newblock \burlalt{\tt [arxiv:2008.12518]}{http://arxiv.org/abs/2008.12518}.

\bibitem[LM10]{Logares.Martens}
M.~Logares and J.~Martens.
\newblock Moduli of parabolic {H}iggs bundles and {A}tiyah algebroids.
\newblock {\em J. Reine Angew. Math.}, 649:89--116, 2010, \doi{10.1515/CRELLE.2010.090}.
\newblock \burlalt{\tt [arxiv:0811.0817]}{http://arxiv.org/abs/0811.0817}.

\bibitem[PMBB23]{MR4556938}
C.~Pauly, J.~Martens, M.~Bolognesi, and T.~Baier.
\newblock The {H}itchin connection in arbitrary characteristic.
\newblock {\em J. Inst. Math. Jussieu}, 22(1):449--492, 2023, \doi{10.1017/S1474748022000196}.
\newblock \burlalt{\tt [arxiv:2002.12288]}{http://arxiv.org/abs/2002.12288}.

\bibitem[Qui85]{MR0783704}
D.~Quillen.
\newblock Determinants of {C}auchy-{R}iemann operators on {R}iemann surfaces.
\newblock {\em Functional Anal. Appl.}, 19(1):31--34, 1985, \doi{10.1007/BF01086022}.
\newblock original: \emph{Funktsional. Anal. i Prilozhen.} 19(1):37--41, 96 (1985).

\bibitem[Sim94]{MR1307297}
C.~T. Simpson.
\newblock Moduli of representations of the fundamental group of a smooth projective variety. {I}.
\newblock {\em Inst. Hautes \'Etudes Sci. Publ. Math.}, (79):47--129, 1994.
\newblock \urlprefix\url{http://www.numdam.org/item?id=PMIHES_1994__79__47_0}.

\bibitem[Sim97]{MR1492538}
C.~Simpson.
\newblock The {H}odge filtration on nonabelian cohomology.
\newblock In {\em Algebraic geometry---{S}anta {C}ruz 1995}, volume 62, Part 2 of {\em Proc. Sympos. Pure Math.}, pages 217--281. Amer. Math. Soc., Providence, RI, 1997, \doi{10.1090/pspum/062.2/1492538}.
\newblock \burlalt{\tt [arxiv:alg-geom/9604005]}{http://arxiv.org/abs/alg-geom/9604005}.

\bibitem[ST04]{MR2106127}
X.~Sun and I.-H. Tsai.
\newblock Hitchin's connection and differential operators with values in the determinant bundle.
\newblock {\em J. Differential Geom.}, 66(2):303--343, 2004.
\newblock \urlprefix\url{http://projecteuclid.org/euclid.jdg/1102538613}.
\newblock \burlalt{\tt [arxiv:math/0309444]}{http://arxiv.org/abs/math/0309444}.

\end{thebibliography}
\end{document}